\documentclass[12pt,a4paper]{article}
\usepackage{amsmath}
\usepackage{amssymb}
\usepackage{amsthm}
\usepackage{amsfonts}
\usepackage{latexsym}

\theoremstyle{plain}
\newtheorem{thm}{Theorem}[section]
\newtheorem{cor}{Corollary}[section]

\newtheorem{prop}{Proposition}[section]
\theoremstyle{definition}

\theoremstyle{remark}

\title{Szeg\"o Kernels and Finite Group Actions}
\author{Roberto Paoletti\footnote{\noindent{\bf Address.}
Dipartimento di Matematica "Ennio De Giorgi", Universit\'a di Lecce,
Via per Arnesano, 73100 Lecce, Italy; {\bf e-mail}:
roberto.paoletti@unile.it }}
\date{}

\begin{document}
\maketitle \footnotesize \noindent {\bf Abstact.} In the context
of almost complex quantization, a natural generalization of
algebro-geometric linear series on a compact symplectic manifold
has been proposed. Here we suppose given a compatible action of a
finite group and consider the linear subseries associated to the
irreducible representations of $G$, give conditions under which
these are base-point free and study properties of the associated
projective morphisms. The results obtained are new even in the complex projective case. \normalsize


\section{Introduction}

Let $(M,\omega)$ be a compact symplectic manifold of dimension $2n$, such that $[\omega ]\in H^2(M,\mathbb{Z})$. Fix $J\in {\cal J}(M,\omega)$ (the contractible space of all almost complex structures on $M$ compatible with $\omega$), and let $h$ and $g={\cal R}(h)$ be the induced hermitian and riemannian structures. There exist an hermitian line bundle $(A,h)$ on $M$
and a unitary covariant derivative $\nabla _A$ on $A$, such that $-2\pi i\omega$ is the curvature of $\nabla _A$. In this set-up, studying the asymptotic spectral properties of a suitable renormalized laplacian, Guillemin and Uribe have introduced privileged spaces of sections $H(M,A^{\otimes k})\subseteq {\cal C}^\infty (M,A^{\otimes k})$; if $J$ is integrable and $A$ holomorphic and ample, these are the spaces of holomorphic sections of $A^{\otimes m}$ \cite{bg}, \cite{gu}, \cite{bu1}. These linear series determine projective embeddings of $M$ enjoying the same metric and symplectic asymptotic properties as in the integrable projective case \cite{bu2}, \cite{z1}, \cite{sz}, \cite{t}.

Suppose that $G$ is a finite group with a symplectic action $\nu :G\times M\rightarrow M$, so that $J$ may be chosen $G$-invariant. Then $\nu$ preserves $g$ and $h$. Assume also that $\nu$ lifts to
a linear action $\tilde \nu:G\times A\rightarrow A$, and that $\tilde \nu$ preserves $h_A$ and $\nabla _A$.
Then $\tilde \nu$ preserves each of the spaces $H(M,A^{\otimes N})$. Let $\rho _i:G\rightarrow {\rm GL}(V_i)$, $1\le i\le c$, be the irreducible representations of $G$; we shall assume that $i=1$ corresponds to the trivial one-dimensional representation. For each $N$, we have a $G$-equivariant decomposition $H(M,A^{\otimes N})=\bigoplus _{i=1}^cH(M,A^{\otimes N})_i$, where $H(M,A^{\otimes N})_i$ consists of a direct sum of copies of $V_i$. It is natural to ask whether the linear series $|H(M,A^{\otimes N})_i|$ are base point free and, if so, what about their asymptotic properties? In this note, we apply  arguments from \cite{bu2} and \cite{z1}, \cite{sz} to these questions.

If $x\in M$, let $G_x=\{g\in G:g\cdot x=x\}$ be its stabiliser. Let $\chi _i :G\rightarrow \mathbb{C}$ be the character of the $i$-th irreducible representation. Let $A_x$ be the fibre of $A$ over $x\in M$. Clearly, $G_x$ acts on $A_x$ and
thus we have a unitary character $\alpha _x:G_x\rightarrow S^1\subset\mathbb{C}^*$. Let $$\gamma _{i,N}(x) :=(\alpha _x^N,\chi _i)_{G_x}=\sum _{g\in G_x}\alpha _x(g)^N\cdot \overline \chi _i(g),\,\,\,\,\,\,\,\,(x\in M,\,1\le i\le c,\, N\in \mathbb{N})$$
$(\,,\,)_{G_x}$ denoting the $L^2$-product with respect to the counting measure on $G_x$.
Note that $\gamma _{i,N}=\gamma _{i,N+|G|}$ for every $i$ and $N$, where $|G|$ denotes the order of $G$. Set
$$B_{i,N}:=\{x\in M:\gamma _{i,N}(x)=0\}=B_{i,N+|G|}\,\,\,\,\,\,\,(1\le i\le c).$$
Clearly, $x\in B_{i,N}$ implies $G_x\neq \{e\}$.

Our first goal is to determine the base locus of the spaces of sections $H(M,A^{\otimes k})_i$ for $k\gg 0$.
In algebro-geometric terminology, the base locus of a vector subspace $W\subseteq {\cal C}^\infty (M,A^{\otimes N})$ is  $${\rm Bs}(|W|):=\{x\in M:\,s(x)=0\,\forall\, s\in W\}.$$

To begin with, we shall prove:

\begin{thm} \label{thm:main1}
Suppose $1\le i\le c$, $0\le r\le |G|-1$, $x\in M$ and $\gamma _{i,r}(x)\neq 0$.
Then for $N\gg 0$, $N\equiv r$ {\rm (mod $|G|$)} there exists a section $s\in H(M,A^{\otimes N})_i$ such that
$s(x)\neq 0$.\end{thm}

This has a number of consequences:

\begin{cor} \label{cor:1} Suppose that the action of $G$ on $M$ is effective. Then $$\dim (H(M,A^{\otimes k})_i)>0$$ for every
$i=1,\ldots,c$ and every $k\gg 0$.
\end{cor}

In fact, it is proved in \cite{pao} that under the same hypothesis
$$\dim (H(M,A^{\otimes k})_i)=\frac{\dim (V_i)^2}{|G|}\cdot \frac{k^n}{n!}\cdot c_1(A)^n+o(k^n).$$

\begin{prop}\label{prop:1} Suppose $1\le i\le c$, $0\le r\le |G|-1$, and $\gamma _{i,r}(x)\neq 0$ for every $x\in M$.
Then $H(M,A^{\otimes k})_i$ globally generates $A^{\otimes k}$ if $k\gg 0$ and $k\equiv r$ {\rm (mod $|G|$)}, that is,
for every $x\in M$ there is $s\in H(M,A^{\otimes k})_i$ such that $s(x)\neq 0$.\end{prop}

\begin{cor} \label{cor:2}
If $k\gg 0$ and $i=1,\ldots,c$, the subspace of $G$-invariant sections $$H(M,A^{\otimes k|G|})^G\subseteq
H(M,A^{\otimes k|G|})$$ globally generates $A^{\otimes k|G|}$.\end{cor}

\begin{cor} \label{cor:3}If $M$ is a complex projective manifold and $A$ is ample,
for every $i=1,\ldots,c$ and $r=0,\ldots,|G|-1$ the base loci ${\rm Bs}\big (|H^0(M,A^{\otimes (r+k|G|)})_i|\big )$
stabilize for $k\gg 0$. Furthermore, for every $k\gg 0$,
$${\rm Bs}\big (|H^0(M,A^{\otimes (r+k|G|)})_i|\big )\subseteq B_{i,r}.$$
\end{cor}

In the reverse direction, it is easily seen that if $G_x=G$ and there exists $s\in {\cal
C}^\infty (M,A^{\otimes N})_i$ with $s(x)\neq 0$, then $$(\alpha
_x^N,\chi _i)_G\neq 0.$$ Therefore,

\begin{cor} \label{cor:4} In the hypothesis of Corollary \ref{cor:3}, suppose in addition that either
$G_x=\{e\}$ or $G_x=G$ for every $x\in G$.
Then $${\rm Bs}\big (|H(M,A^{\otimes N})_i|\big )= B_{i,N}$$ for $i=1,\ldots,c$ and $N\gg 0$. \end{cor}

In the almost complex case, for any $i=1,\ldots,c$ and $r=0,\ldots,|G|-1$
we may still define the $(i,r)$-th equivariant {\rm asymptotic} base locus of
$A$ as
\begin{eqnarray*}{\rm Bs}(A,i,r)_\infty =:\left  \{x\in M:\, \forall s>0\, \exists \,k>s,\, k\equiv r\mbox{ (mod.$|G|$)}\right .\\\mbox{ such that }
\left . x\in {\rm Bs}\big (|H(M,A^{\otimes k})_i|\big )\right \}.\end{eqnarray*}
The general case (symplectic, almost complex) of Corollary \ref{cor:3} is then

\begin{cor} \label{cor:5} In the above situation,
$${\rm Bs}(A,i,r)_\infty \subseteq B_{i,r}.$$
If furthermore $K\subset M$ is any compact subset with $K\cap B_{i,r}=\emptyset$,
then $$K\cap {\rm Bs}\big (|H(M,A^{\otimes k})_i|\big )=\emptyset$$ for all $k\gg 0$ with $k\equiv r\,\,{\rm (mod.\, |G|)}$.
\end{cor}

Next, if $B_{i,r}=\emptyset$ there are associated projective morphisms
$$\Phi _{i,r+k|G|}:M\rightarrow \mathbb{P}(H(M,A^{\otimes (r+k|G|)})^*_i),$$
and we now consider their
asymptotic properties as $k\rightarrow +\infty$.

\begin{thm} \label{thm:2}
Suppose $B_{i,r}=\emptyset$ for some $1\le i\le c$ and $0\le r\le |G|-1$.
Let $U\subseteq M$ be the open subset of $M$ where the order $|G_x|$ is locally constant. Suppose
$U'\subset U$ is open with $\overline {U'}\subset U$. Then $\Phi _{i,r+k|G|}$ is an immersion on $U'$ for $k\gg 0$. \end{thm}

\begin{cor} $|H(M,A^{\otimes N})^G|$ is base point free and $\Phi _{1,N}$ is an
immersion on compact subsets of $U$ if
$N\gg 0$ and $\sum _{g\in G_x}\alpha _x(g)^N\neq 0$ for every $x\in G$.\end{cor}

In general $\Phi _{i,N}$ is not injective; for example it is constant on every
orbit for any $G$ if $i$ corresponds to the trivial representation, or for any $i$ if $G$ is abelian. We may still ask,
however, if in these cases
points in different orbits have different images under $\Phi _{i,N}$.

Let $d_G:M\times M\rightarrow \mathbb{R}$ be the {\it orbit distance}:
$$d_G(x,y):=\min \{d(gx,y):g\in G\} \,\,\,\,\,\,\,\,\mbox{($x,y\in M$)}.$$ Clearly, $d_G(x,y)>0$ if and only
if $x\not \in G\cdot y$.

\begin{prop} \label{prop:main2} Assume that either $G$ is abelian, or $G$ is arbitrary and $i=1$.
Let $U\subseteq M$ be as in Theorem \ref{thm:main1}, $N\in \mathbb{N}$ and suppose that $B_{i,N}=\emptyset$ and that $\gamma _{i,N}$ is constant on $W$. Let $K\subseteq W$ be a compact subset.
There exists $k_0\in \mathbb{N}$ such that if $k\ge k_0$, $x,y\in K$ and $d_G(x,y)>0$ then
$$\Phi _{i,N+k|G|}(x)\neq \Phi _{i,N+k|G|}(y).$$
\end{prop}

\begin{cor} If the action of $G$ on $M$ is free, then $\Phi_{i,N}$ is well-defined and
an embedding $M/G\hookrightarrow \mathbb{P}\big ({H(M,A^{\otimes
k})^G}^*\big )$ for any $i=1,\ldots ,c$ and $N\gg 0$.
\end{cor}

Similar statements hold for the asymptotic metric and almost
complex properties, in the vein of theorem 1.1 of \cite{bu2}.

\section{Proofs}
{\it Proof of Theorem \ref{thm:main1}.}
We recall some notation from \cite{bu2}, \cite{z1}, \cite{sz}.
Let $A^*=A^{-1}$ be the dual line bundle with the induced hermitian
stucture $h_{A^*}$, and let $A^*\supset \mathbb{S}\stackrel{\pi}{\rightarrow}M$
be the unit circle bundle, a strictly pseudoconvex domain. Given the connection,
$\mathbb{S}$ has natural riemannian and almost CR structures. We shall identify
functions and half-forms throughout.

As $\mathbb{S}$ is a principal $S^1$-bundle, ${\cal C}^\infty
(\mathbb{S})=\bigoplus _{N\in \mathbb{Z}}{\cal C}^\infty
(\mathbb{S})_N$, where ${\cal C}^\infty (\mathbb{S})_N$ is the
$N$-th isotype for the $S^1$-action. In the standard manner, we
shall identify ${\cal C}^\infty (M,A^{\otimes N})$ and ${\cal
C}^\infty (\mathbb{S})_N$. Set $H(\mathbb{S}):=\bigoplus _{N\in
\mathbb{N}}H(\mathbb{S})_N$, where $H(\mathbb{S})_N\cong
H(M,A^{\otimes N})$ under this identification; in the integrable
projective case, $H(\mathbb{S})$ is the Hardy space of boundary
values of holomorphic functions on $A^*$. Let $\Pi
:L^2(\mathbb{S})\rightarrow H(\mathbb{S})$ be the orthogonal
projector and $\tilde \Pi \in {\cal D}'(\mathbb{S}\times
\mathbb{S})$ its Schwartz kernel; decompose it as $\tilde \Pi
=\bigoplus _{N\in \mathbb{N}}\tilde \Pi _N$, where $\tilde \Pi
_N\in {\cal C}^\infty (\mathbb{S}\times \mathbb{S})$ is the $N$-th
Fourier coefficent. We have $\tilde \Pi _N(x,y)=\sum
_{i=0}^{d_N}s_i^N(x)\otimes \overline s_i^N(y)$, where
$\{s_0^N,\ldots,s^N_{d_N}\}$ is an orthonormal basis of
$H(\mathbb{S})_N$. Let $\tilde \Phi _{i,N}:\mathbb{S}\rightarrow
H(M,A^{\otimes N})^*$ be the coherent state map, given by
evaluation, which is a lifting of $\Phi _{i,N}$ when the latter is
defined. Then $\tilde \Pi_N(p,q)=\big (\tilde \Phi
_{i,N}(p),\tilde \Phi _{i,N}(q)\big )$ ($p,q\in \mathbb{S}$),
where $(\cdot \,,\cdot)$ denotes the $L^2$-hermitian product on
$H(M,A^{\otimes N})^*$.

The induced action of $G$ on $A^*$ preserves $\mathbb{S}$ and the
riemannian and almost CR structures on $\mathbb{S}$, and the
isomorphisms $H(\mathbb{S})_N\cong H(M,A^{\otimes N})$ are
$G$-equivariant. For $N\gg 0$, we have $G$-equivariant
decompositions $H(\mathbb{S})_N= \bigoplus _iH(\mathbb{S})_{i,N}$,
where $H(\mathbb{S})_{i,N}$ is the factor consisting of a direct
sum of copies of $V_i$, $1\le i\le c$. Similarly,
$H(\mathbb{S})=\bigoplus _iH(\mathbb{S})_i$. We shall implicitly
identify $H(\mathbb{S})_N$ and $H(\mathbb{S})_{i,N}$ with
$H(M,A^{\otimes N})$ and $H(M,A^{\otimes N})_i$, respectively. For
each $i$, let $\Pi _i:L^2(\mathbb{S})\rightarrow H(\mathbb{S})_i$
denote the orthogonal projection and let $\tilde \Pi _i\in {\cal
D}'(\mathbb{S}\times \mathbb{S})$ be its Schwartz kernel. For each
$i$ and $N$, let $\Pi _{i,N}:L^2(\mathbb{S})\rightarrow
H(\mathbb{S})_{i,N}$ be the orthogonal pojection and $\tilde \Pi
_{i,N}$ its Schwartz kernel, the $N$-th Fourier coefficient of
$\tilde \Pi _i$: if $\{s_0^{(i,N)},\ldots,s_{d_{i,N}}^{(i,N)}\}$
is an orthonormal basis of $H(\mathbb{S})_{i,N}$, then $$\tilde \Pi
_{i,N}(p,q)=\sum _{j=0}^{d_{i,N}}s_j^{(i,N)}(p)\otimes \overline
{s_j^{(i,N)}(q)} \,\,\,\,\,\,\,\,\mbox{($p,q\in \mathbb{S}$)}.$$
Clearly, $\tilde \Pi =\sum _{i=1}^c\tilde \Pi _i$.
Notice that the Fourier components of the total and equivariant Szeg\"o kernels,
$\Pi _N$ and $\Pi _{i,N}$, when restricted to the diagonal in $\mathbb{S}\times \mathbb{S}$
descend to well-defined smooth functions on
the diagonal in $M\times M$, that is, we may write with some abuse of language $\Pi _N(p,p)=\Pi _N(x,x)$
and $\Pi _{i,N}(p,p)=\Pi _{i,N}(x,x)$ for any $p\in \mathbb{S}$ and $x\in M$ with
$\pi (p)=x$. This will be done implicitly below.

By the projection formula, for each $i=1,\ldots,c$ we have
\begin{eqnarray*}\tilde \Pi _{i,N}=\sum
_{j=0}^{d_N}\Pi _i (s_j^N)\otimes \overline s_j^N=(\dim
(V_i)/|G|)\cdot \sum _g\sum _j\overline \chi _i(g)\rho
(g)(s_j^N)\otimes \overline s_j^N,\end{eqnarray*} where $\rho :G\rightarrow {\rm
GL}(H(\mathbb{S})_N)$ is the induced representation; explicitly,
$\rho (g)\sigma =\sigma \circ g^{-1}$ ($g\in G,\, \sigma\in
H(\mathbb{S})_N$), where we view $g^{-1}$ as a contactomorphism of
$\mathbb{S}$. Thus,
\begin{eqnarray*}\label{eqn:0}\tilde \Pi _{i,N}(p,q)=(\dim (V_i)/|G|)\cdot \sum
_g\sum _j\overline \chi _i(g)s_j^N(g^{-1}p) \overline
s_j^N(q)\\=(\dim (V_i)/|G|)\cdot \sum _g\overline \chi _i(g)\tilde
\Pi _N(g^{-1}p,q).\end{eqnarray*} On the diagonal, $\tilde \Pi
_{i,N}(p,p)=(\dim (V_i)/|G|)\cdot \sum _g\overline \chi
_i(g)\tilde \Pi _N(g^{-1}p,p)$. Let $d$ be the geodesic distance
function on $M$ and also its pull-back $d\circ \pi$  to
$\mathbb{S}$. If $x\in M$ and $G\cdot x\neq \{x\}$, set $a_x:=\min
\{d(gx,x):g\in G\setminus G_x\}$. Suppose $p\in \mathbb{S}$,
$x=\pi (p)$. Then
\begin{eqnarray*}\tilde \Pi _{i,N}(p,p)=(\dim (V_i)/|G|)\cdot \sum _{g\in G_x}\overline
\chi _i(g)\tilde \Pi _N(g^{-1}p,p)+\\(\dim (V_i)/|G|)\cdot \sum
_{g\not \in G_x}\overline \chi _i(g)\tilde \Pi
_N(g^{-1}p,p).\end{eqnarray*} By Lemma 4.5 of \cite{bu2}, the
latter term is bounded in absolute value by $C\tilde \Pi
_N(p,p)e^{-a_x^2 N/2}+O(N^{(n-1)/2}),$ where $C$ is independent of
$x$ and $N$. By (13) of \cite{sz0} and the definition of dual
action, $\tilde \Pi _N(g^{-1}p,p)=\alpha _x(g)^N\tilde \Pi
_N(p,p)$ if $g\in G_x$. Thus the former term is
\begin{eqnarray*}(\dim (V_i)/|G|)\cdot \Big [\sum _{g\in G_x}
\overline \chi _i(g)\alpha _x(g)^N\Big ]\tilde \Pi _N(p,p)= \\(\dim (V_i)/|G|)
\cdot (\alpha _x^N,\chi _i)_{G_x}\cdot \tilde \Pi _N(p,p).\end{eqnarray*}
Given the asymptotic expansion for $\tilde \Pi _N(p,p)$ in \cite{bu2} and \cite{z1},
$\tilde \Pi _{i,N}(p,p)\neq 0$ if $N\gg 0$, $x\not \in B_{i,N}$.
This clearly implies the statement.

\bigskip

\noindent {\it Proof Of Corollary \ref{cor:1}.} Let $V\subseteq M$ be the locus of points
with non-trivial stabilizer. By Theorem 8.1 on page 213 of \cite{s} and because the action
is effective, $V$ is a union of proper submanifolds of $M$. If $x\in M\setminus V$, then
$G_x=\{e\}$ and therefore $\gamma _{i,k}(x)=\dim (V_i)\neq 0$ for every $i$ and $N$. By the
Theorem, there exists $s\in H(M,A^{\otimes k})_i$ with $s(x)\neq 0$ if $k\gg 0$.

\bigskip

Before coming to the proof of Proposition \ref{prop:1}, let us dwell on the previous
descrition of the equivariant Szeg\"o kernels $\tilde \Pi _{i,k}$ restricted to the diagonal.
As is well-known, the wave front of the Szeg\"o kernel $\Pi$ is
$$\Sigma =\left \{\left ((p,p),(r\alpha _p,-r\alpha _p)\big ):p\in \mathbb{S},r>0\right )\right \}\subseteq T^*\left (\mathbb{S}\times \mathbb{S}\right )\setminus \{0\}.$$
In the notation of \cite{bg}, \cite{bu2} we have in fact $\Pi \in J^{1/2}(\mathbb{S}\times \mathbb{S},\Sigma)$.
Now we have seen that
$$\tilde \Pi
_{i,N}(p,p)=(\dim (V_i)/|G|)\cdot \sum _{g\in G}\overline \chi
_i(g)\tilde \Pi _N(g^{-1}p,p).$$
For any $g\in G$ let $\alpha _g: \mathbb{S}\times \mathbb{S}\rightarrow \mathbb{S}\times \mathbb{S}$
be the diffeomorphism $(p,q)\mapsto (g\,p,q)$, and let $\tilde \Pi _g=\tilde \Pi \circ \alpha _g^*\in {\cal D}'(\mathbb{S}\times \mathbb{S})$, where
$\alpha _g^*$ denotes pull-back of functions under $\alpha _g$. Then
$\tilde \Pi _g\in J^{1/2}\big (\mathbb{S}\times \mathbb{S},\alpha _g^*(\Sigma)\big )$ and $\tilde \Pi _k(g\, p,q)$ is the
$k$-th Fourier component of $\tilde \Pi _g$, for every integer $k$. One can then see, arguing as in the proofs of Lemmas
3.5 and 3.6 of \cite{bu2}, that $k^{-n}\tilde \Pi _k(g\, p,p)$ is bounded in ${\cal C}^1$ norm, say, for every $g\in G$.
The same then holds for $k^{-n}\tilde \Pi _{i,k}(x,x)$.

\bigskip

\noindent {\it Proof of Proposition \ref{prop:1}.} By the above, in the hypothesis of the Proposition
for every $x\in M$ there exists
$k_x\in \mathbb{N}$ such that $x\not\in {\rm Bs}\big (|H(M,A^{\otimes k})_i|\big )$ for every $k\ge k_x$.
We now make the stronger claim that for every $x\in M$ there exist an open neighbourhood
$U_x$ of $x$ and $k_x\in \mathbb{N}$ such that $U_x\cap {\rm Bs}\big (|H(M,A^{\otimes k})_i|\big )=\emptyset$
for every $k\ge k_x$. The statement will follow given the compactness of $M$.

If the claim was false, there would exist $x\in M$ and sequences $k_j\in \mathbb{N}$ and $x_j\in M$
($j=1,2,\ldots$) with $k_j\equiv r$ (mod $|G|$), $k_j\rightarrow +\infty$ and $x_j\rightarrow x$, such that $x_j\in {\rm Bs}\big (|H(M,A^{\otimes k_j})_i|\big )$ for every $j$. Thus,
$$\tilde \Pi _{i,k_j}(x_j,x_j)=0\,\,\,\,\,\,\,\mbox{($j=1,2,\ldots$)}$$
while $$\tilde \Pi _{i,k_j}(x,x)=\frac{\dim (V_i)}{|G|}\cdot \gamma _{i,r}(x)\cdot \tilde \Pi _{k_j}(x,x)+{\rm L.O.T.},$$
where L.O.T. denotes lower order terms in $k_j$. Thus, $k_j^{-n}\, \tilde \Pi _{i,k_j}(x,x)$ is bounded away from zero
and therefore the derivatives in $x$ of the sequence of functions $k_j^{-n}\, \tilde \Pi _{i,k_j}(x',x')$ are unbounded, a contradiction.

\bigskip

\noindent
{\it Proof of Corollary \ref{cor:2}.} Let us agree that the irreducible representation corresponding to $i=1$ is just the trivial representation, so that $$H(M,A^{\otimes N})_1=H(M,A^{\otimes N})^G$$ for every integer $N$. Then $\overline \chi _1(g)=1$ for every $g\in G$. Furthermore, for every $x\in M$, $g\in G_x$ and $k\in \mathbb{N}$ we have $\alpha _x^{k|G|}(g)=1$. Thus $$\gamma _{1,k|G|}(x)=|G_x|\neq 0\,\,\,\mbox{for every $x\in M$,}$$ and the statement follows from  Proposition \ref{prop:1}.

\bigskip

\noindent {\it Proof of Corollary \ref{cor:3}.} If $M$ is a complex projective manifold and $A$ is ample,
we have section multiplication maps
$$H^0(M,A^{\otimes \ell})^G\otimes H^0(M,A^{\otimes m})_i\longrightarrow H^0(M,A^{\otimes (\ell +m)})_i$$
for every $i=1,\ldots,c$ and integers $\ell,\, m$. Thus, for any residue class $0\le r\le |G|-1$ and any sequence
of positive integers
$k_i\gg 0$, by
Corollary \ref{cor:2} we have a descending chain of base loci:
\begin{eqnarray*}{\rm Bs}\big (|H^0(M,A^{\otimes r})_i|\big )\supseteq {\rm Bs}\big (|H^0(M,A^{\otimes (r+k_1|G|)})_i|\big )\supseteq \\
{\rm Bs}\big (|H^0(M,A^{\otimes (r+(k_1+k_2)|G|)})_i|\big )\supseteq \ldots.\end{eqnarray*}
This implies the first statement. The rest is obvious.

\bigskip

\noindent {\it Proof of Corollary \ref{cor:4}.}
If $G_x=G$ and $k\equiv r$ (mod. $|G|$), then
$$\tilde \Pi _{i,k}(x,x)=\frac{\dim (V_i)}{|G|}\cdot \gamma _{i,r}(x)\cdot \tilde \Pi _k(p,p).$$
Thus, if $\gamma _{i,r}(x)=0$ then $s(x)=0$ for every $s\in H(M,A^{\otimes k})_i$.

\bigskip

\noindent {\it Proof of Corollary \ref{cor:5}.} The first statement follows from Theorem \ref{thm:main1},
while the second is a consequence of the proof of Proposition \ref{prop:1}.

\bigskip

\noindent {\it Proof of Theorem \ref{thm:2}.}
Suppose
$B_{i,N}=\emptyset$ so that, perhaps after replacing $N$ by
$N+k|G|$ for $k\gg 0$, $|H(\mathbb{S})_{i,N}|$ is base point free.
The claim is that if $U'\subset U$ is open with compact closure in
$U$ and $N\gg 0$, then $\Phi _{i,N}$ is an immersion on $U'$. We
shall be done by proving that $N^{-1}\Phi _{i,N}^*(\omega _{\rm
FS}^{(N)})-\omega =O(1/N)$ on connected compact subsets of $U$,
where $\omega _{\rm FS}^{(N)}$ is the Fubini-Study symplectic form
on $\mathbb{P}\big (H(M,A^{\otimes k})^*)$. In turn, this will
follow if we prove that $N^{-1}\tilde \Phi _{i,N}^*(\tilde \omega
_N)-\pi ^*(\omega)=O(1/N)$ on horizontal vectors, over compact
subsets of $\mathbb{S}$; here $\tilde \omega _N=\frac i2\overline
\partial
\partial \log |\xi|^2$ on $H(M,A^{\otimes k})^*\setminus \{0\}$
(with its natural hermitian structure), and $\pi
:\mathbb{S}\rightarrow M$ is the projection.

Now, if $d^1$ and $d^2$ denote exterior differentiation on the
first and second component of $\mathbb{S}\times \mathbb{S}$,
respectively, then $N^{-1}\tilde \Phi _{i,N}^*\tilde \omega
_N={\rm diag}^*(d^1d^2\log \tilde \Pi _{i,N})$, where ${\rm
diag}:\mathbb{S}\rightarrow \mathbb{S}\times \mathbb{S}$ is the
diagonal map (\cite{sz}, proof of theorem 3.1 (b)). If $x,y\in M$
lie in the same connected component $V$ of $U$, $G_y=G_x$. Thus
$b_x:=(\alpha _x^N,\chi _i)_{G_x}$ is constant on $V$, say equal
to $b_V$. Hence, if $p,q\in \pi ^{-1}(V)$ and $x=\pi (p)$,
\begin{equation}\label{eqn:baba'}
\tilde \Pi _{i,N}(p,q)=\frac {\dim (V_i)}{|G|}\cdot \Big
\{b_V\cdot \tilde \Pi _N(p,q)+ \sum _{g\not \in G_x}\overline \chi
_i(g)\tilde \Pi _N(gp,q)\Big \}.\end{equation} By the proof of
theorem 3.1 (b) of \cite{sz}, $(i/2N){\rm diag}^*\big (d^1d^2\log
\Pi _N\big )\rightarrow \pi ^*\omega$ in ${\cal C}^k$-norm for any
$k$ on $M$. Therefore, we shall be done by proving that
\begin{equation}\label{eqn:pepe'}N^{-1}d_1d_2
\big (\tilde \Pi _N(gp,q)/\tilde \Pi _N(p,q)\big
)\rightarrow 0\end{equation} and \begin{equation}N^{-1}d_1(\big
(\tilde \Pi _N(gp,q)/\tilde \Pi _N(p,q)\big )\wedge d_2(\big
(\tilde \Pi _N(g'p,q)/\tilde \Pi _N(p,q)\big )\rightarrow
0\label{eqn:parapa'}\end{equation} for $g,g'\not \in G_x$ near
compact subsets of ${\rm diag}(V)$.

Let then $K\subset V$ be a compact subset, and suppose $x\in K$,
$g\not \in G_x$, and $u,v\in T_xM$ have unit length. Let $U,V$ be
horizontal vector fields of unit length on $\mathbb{S}$, of unit
length near $\mathbb{S}_x$ and extending the horizontal lifts of
$u$ and $v$. We want to estimate $N^{-1}U_1\circ V_2(\tilde \Pi
_N(gp,q)/\tilde \Pi _N(p,q))$ over $K$, where $U_1=(U,0)$ and
$V_2=(0,V)$ are horizontal vector fields on $\mathbb{S}\times
\mathbb{S}$.

Let us consider again the distribution $\tilde \Pi _g=
\alpha _g ^*\tilde \Pi \in J^{1/2}(\mathbb{S}\times
\mathbb{S},g^*\Sigma)$,
discussed before the proof of Proposition \ref{prop:1}.
If $P$ is a horizontal differential operator of degree $\ell$ on
$\mathbb{S}\times \mathbb{S}$, its principal symbol vanishes on
$g^*\Sigma$ and therefore $P(\tilde \Pi _g)\in
J^{(\ell+1)/2}(\mathbb{S}\times \mathbb{S},\alpha _g^*\Sigma)$. As in
\cite{bu2} Lemma 4.5, for $k\in \mathbb{N}$ we can find $\nu
_{{}_{g,P,k}}\in {\cal C}^\infty (\mathbb{S})$, having an
asymptotic development $\nu _{{}_{g,P,k}}(p)=\sum _{j=0}^\infty
k^{n+(\ell -j)/2}f_{{}_{g,P,k}}^{(j)}(p)$, and real phase
functions $\alpha _{{}_{g,P,k}}\in {\cal C}^\infty
(\mathbb{S}\times \mathbb{S})$ such that $$G(p,q)=\sum _k\nu
_{{}_{g,P,k}}(p)e^{i \alpha
_{{}_{g,P,k}}(p,q)}e^{-kd(gp,q)^2/2}\in J^{(\ell
+1)/2}(\mathbb{S}\times \mathbb{S},\alpha _g^*\,\Sigma)$$ and $P(\tilde \Pi _g)-G\in
J^{\ell /2}(\mathbb{S}\times \mathbb{S},\alpha _g^*\, \Sigma).$ Since
$P(\tilde \Pi _g)$ has definite (even) parity, we may assume
without loss that so does $G$. Therefore, the above asymptotic
expansions may be assumed to go down by integer steps: $\nu
_{{}_{g,P,k}}(p)=\sum _{j=0}^\infty k^{n+\ell/2
-j}f_{{}_{g,P,k}}^{(j)}(p)$, and \begin{equation}
\label{eqn:lola}\big |P\big (\tilde \Pi _N(gp,q)\big )\big |=\nu
_{{}_{g,P,0}}(p)\cdot
e^{-Nd(gp,q)^2/2}+O(N^{n+\ell/2-1}).\end{equation}

Because $K\subset U$ is compact and $g\not \in G_x$ for $x\in K$,
there is $\epsilon
>0$ such that $d(gp,p)>\epsilon$ for all $p\in \pi ^{-1}(K)$.
Thus, $P(\tilde \Pi _N^{(g)})(p,p)=O(N^{m+(\ell -1)/2})$ on $\pi
^{-1}(K)$. Developping $N^{-1}U_1\circ V_2(\tilde \Pi
_N(gp,q)/\tilde \Pi _N(p,q))$, we see that $N^{-1}U_1\circ
V_2(\tilde \Pi _N(gp,q)/\tilde \Pi _N(p,q))=O(1/N)$ over $K$,
uniformly in $U$ and $V$. This proves (\ref{eqn:pepe'}); the proof
of the other estimate is similar.

\bigskip

\noindent {\it Proof of Proposition \ref{prop:main2}.}
Notation being as above, we may assume that
$K$ is $G$-invariant. Suppose then, by contradiction, that for a
sequence $k_j\rightarrow +\infty$ we can find $x_{k_j},y_{k_j}\in
K$ with $d_G(x_{k_j},y_{k_j})>0$ and $\Phi
_{i,N+k_j|G|}(x_{k_j})=\Phi _{i,N+k_j|G|}(y_{k_j})$. Set
$N_j=N+k_j|G|$.

I claim that $d_G(x_{k_j},y_{k_j})\le C/\sqrt {N_j}$. Following
\cite{bu2}, proof of Corollary 4.6, pick $p_{k_j}\in \pi
^{-1}(x_{k_j})$, $q_{k_j}\in \pi ^{-1}(y_{k_j})$. Then $\tilde
\Phi _{i,N_j}(x_{k_j})=\lambda _j\tilde \Phi _{i,N_j}(y_{k_j})$
for some $\lambda _j\in \mathbb{C}$; it follows that $||\tilde
\Phi _{i,N_j}(p_{k_j})||^2=|\lambda _j|^2\cdot ||\tilde \Phi
_{i,N_j}(q_{k_j})||^2$. However, $||\tilde \Phi
_{i,N_j}(p)||^2=\tilde \Pi _{i,N_j}(p,p)$ ($p\in \mathbb{S}$), and
therefore by (\ref{eqn:baba'}) above $|\lambda
_j|=1+O(N_j^{-1/2})$. We also have $|\lambda _j|\tilde \Pi
_{i,N_j}(p_{k_j},p_{k_j})=|\tilde \Pi _{i,N_j}(p_{k_j},q_{k_j})|$,
and on the other hand, again by (\ref{eqn:baba'}), $$|\tilde \Pi
_{i,N_j}(p_{k_j},q_{k_j})|\le C |\tilde \Pi
_{i,N_j}(p_{k_j},p_{k_j})|e^{-N_jd_G(p,q)^2/2}+O(k_j^{n-1/2}).$$
We conclude that $d_G(p_{k_j},q_{k_j})\le C/\sqrt {k_j}$, as
claimed. Hence, after replacing $x_{k_j}$ by $g_j\cdot x_{k_j}$
for a suitable $g_j\in G$, we may assume $d(x_{k_j},y_{k_j})\le
C/\sqrt {N_j}$ and $d(x_{k_j},y_{k_j})=d_G(x_{k_j},y_{k_j})$ for
every $j$.

Since $d(gx,x)>\epsilon$ for some fixed $\epsilon >0$ and all
$x\in K$ and $g\not \in G_x$, $x_{k_j}$ is the only point in
$G\cdot x_{k_j}$ minimizing the distance from $y_{k_j}$, for every
$j$.

We may now apply the argument of the proof of theorem 3.2 (b) of \cite{sz}, with minor changes.

\end{document}